\documentclass[9pt]{article}

\usepackage[cmex10]{amsmath}
\usepackage{amssymb}
\usepackage{graphicx}
\usepackage[numbers,sort&compress]{natbib}
\usepackage[squaren]{SIunits}
\usepackage{needspace}
\usepackage{pmat}

\usepackage{tikz}
\usetikzlibrary{shadows}
\usetikzlibrary{shapes}
\usetikzlibrary{decorations}
\usepackage{framed}
\usepackage[framed]{ntheorem}

\theoremclass{Theorem}
\theoremstyle{plain}

\newshadedtheorem{theorem}[Theorem]{Theorem}

\renewcommand{\d}[1]{\;\mathrm{d}#1}
\renewcommand{\t}[1]{\mathrm{T}#1}

\newcommand{\erf}{\operatorname{erf}}

\newcommand{\N}{\mathcal{N}}

\newcommand{\CN}{\mathcal{CN}}
\newcommand{\V}{\mathcal{V}}
\newcommand{\I}{\mathcal{I}}

\newcommand{\R}{\mathcal{R}}
\newcommand{\Ncdf}{\mathcal{N}\hspace{-0.1cm}\operatorname{cdf}}

\renewcommand{\d}[1]{\;\mathrm{d}#1}
\renewcommand{\t}[1]{\mathrm{T}#1}
\newcommand{\h}[1]{\mathrm{H}#1}

\title{Calculation of Some Expected Values for Parameterized Mean Model with Gaussian Noise}
\author{Umut Orguner\thanks{Department of Electrical \& Electronics Engineering, Middle East Technical University, Ankara, Turkey}}
\begin{document}
\setlength{\interdisplaylinepenalty}{2500}
\maketitle
Consider the measurement model
\begin{align}
y=g(x)+v\label{eqn:commongeneralexample}
\end{align}
where
\begin{itemize}
\item $x\in\mathbb{R}$ is the unknown scalar we would like to estimate described by the prior distribution $\N(x;0;\sigma_x^2)$
where the notation $\N(\cdot;\bar{x},\Sigma_x)$ denotes a (real-variate) Gaussian density with mean $\bar{x}$ and covariance $\Sigma_x$.
\item $y\in\mathbb{R}^{n_y}$ is, in general, a complex measurement vector;
\item $g(\cdot): \mathbb{R}\rightarrow \mathbb{C}^{n_y}$ is, in general, a complex-valued observation function;
\item $v\in\mathbb{C}^{n_y}$ is circular symmetric complex Gaussian measurement noise with zero-mean  and covariance $\sigma_v^2 I_{n_y}$ where $I_{n_y}$ denotes an identity matrix of size $n_y\times n_y$. The noise $v$ is assumed independent of $x$.
\end{itemize}

In this document we are going to derive analytical formulae for the following functions.
\begin{align}
\mu_{y,x}(s_1,s_2,h_1,h_2)\triangleq& E\Big[\min(L^{s_1}(y,x+h_1,x),1)\min(L^{s_2}(y,x+h_2,x),1)\Big],\\
\mu_{y|x}(s_1,s_2,h_1,h_2,x)\triangleq& E_y\Big[\min(L^{s_1}_1(y,x+h_1,x),1)\min(L^{s_2}_1(y,x+h_2,x),1)\Big|x\Big],\\
\mu_{x}(s_1,s_2,h_1,h_2)\triangleq& E_x\Big[\min(L^{s_1}_2(x+h_1,x),1)\min(L^{s_2}_2(x+h_2,x),1)\Big].
\end{align}
where $s_1,s_2\in\mathbb{Z}$, $h_1,h_2\in\mathbb{R}^{n_x}$ and
\begin{align}
L(y,x,\xi)\triangleq& \frac{p_{\tilde{y},\tilde{x}}(y,x)}{p_{\tilde{y},\tilde{x}}(y,\xi)},\qquad L_1(y,x,\xi)\triangleq \frac{p_{\tilde{y}|\tilde{x}}(y|x)}{p_{\tilde{y}|\tilde{x}}(y|\xi)},\qquad L_2(x,\xi)\triangleq\frac{p_{\tilde{x}}(x)}{p_{\tilde{x}}(\xi)}.
\end{align}
\section{Calculation of $\mu_{y,x}(s_1,s_2,h_1,h_2)$}
\label{sec:muyx}
We first calculate $L(y,x+h,x)$ as
\begin{align}
L(y,x+h,x)=&\frac{\CN(y;g(x+h),\sigma_v^2I_{n_y})\N(x+h;0;\sigma_x^2)}{\CN(y;g(x),\sigma_v^2I_{n_y})\N(x;0;\sigma_x^2)}\\
=&\exp\bigg[\frac{2}{\sigma_v^2}\R\left\{v^\h d(x,h)\right\}-b(x,h)\bigg]
\end{align}
where the notation $\CN(y;\bar{y},\Sigma_y)$ denotes a circular symmetric complex Gaussian density with mean $\bar{y}$ and covariance $\Sigma_y$ and
\begin{align}
d(x,h)\triangleq& g(x+h)-g(x),\label{eqn:defintiondxh}\\
b(x,h)\triangleq& \frac{1}{\sigma_v^2}\|d(x,h)\|^2+\frac{1}{\sigma_x^2}x^\t h+\frac{1}{2\sigma_x^2}\|h\|^2.\label{eqn:defintionbxh}
\end{align}
Then we have
\begin{align}
\mu_{y,x}&(s_1,s_2,h_1,h_2)\nonumber\\
=&\int\int \min\bigg(\exp\bigg[\frac{2s_1}{\sigma_v^2}\R\left\{v^\h d(x,h_1)\right\}-s_1b(x,h_1)\bigg],1\bigg)\nonumber\\
&\times\min\bigg(\exp\bigg[\frac{2s_2}{\sigma_v^2}\R\left\{v^\h d(x,h_2)\right\}-s_2b(x,h_2)\bigg],1\bigg)p(v)\d v p(x)\d x.\label{eqn:muxyderivation1}
\end{align}
We are first going to handle the inner integral on the right hand side of~\eqref{eqn:muxyderivation1} which we call as $\I_1$ as follows.
\begin{align}
\I_1&\triangleq\int\min\bigg(\exp\bigg[\frac{2s_1}{\sigma_v^2}\R\big\{v^\h d(x,h_1)\big\}-s_1b(x,h_1)\bigg],1\bigg)\nonumber\\
&\times\min\bigg(\exp\bigg[\frac{2s_2}{\sigma_v^2}\R\big\{v^\h d(x,h_2)\big\}-s_2b(x,h_2)\bigg],1\bigg)p(v)\d v\nonumber\\
=&\int_{\V_1}p(v)\d v\nonumber\\
&+\int_{\V_2} \exp\bigg[\frac{2s_2}{\sigma_v^2}\R\big\{v^\h d(x,h_2)\big\}-s_2b(x,h_2)\bigg]p(v)\d v\nonumber\\
&+\int_{\V_3} \exp\bigg[\frac{2s_1}{\sigma_v^2}\R\big\{v^\h d(x,h_1)\big\}-s_1b(x,h_1)\bigg]p(v)\d v\nonumber\\
&+\int_{\V_4}
\exp\bigg[\frac{2s_1}{\sigma_v^2}\R\big\{v^\h d(x,h_1)\big\}-s_1b(x,h_1)\bigg]\exp\bigg[\frac{2s_2}{\sigma_v^2}\R\big\{v^\h d(x,h_2)\big\}-s_2b(x,h_2)\bigg]p(v)\d v\label{eqn:muxyderivation2}
\end{align}
where the sets $\V_1$, $\V_2$, $\V_3$ and $\V_4$ are defined as follows.
\begin{subequations}
\begin{align}
\V_1&\triangleq\Big\{v\in\mathbb{C}^{n_y}\Big|{\scriptstyle \R\big\{v^\h d(x,h_1)\big\}\ge\frac{\sigma_v^2 b(x,h_1)}{2} \& \R\big\{v^\h d(x,h_2)\big\}\ge\frac{\sigma_v^2b(x,h_2)}{2}} \Big\},\\
\V_2&\triangleq\Big\{v\in\mathbb{C}^{n_y}\Big|{\scriptstyle \R\big\{v^\h d(x,h_1)\big\}\ge\frac{\sigma_v^2 b(x,h_1)}{2} \& \R\big\{v^\h d(x,h_2)\big\}< \frac{\sigma_v^2b(x,h_2)}{2}} \Big\},\\
\V_3&\triangleq\Big\{v\in\mathbb{C}^{n_y}\Big|{\scriptstyle \R\big\{v^\h d(x,h_1)\big\}< \frac{\sigma_v^2 b(x,h_1)}{2} \& \R\big\{v^\h d(x,h_2)\big\}\ge\frac{\sigma_v^2b(x,h_2)}{2}} \Big\},\\
\V_4&\triangleq\Big\{v\in\mathbb{C}^{n_y}\Big|{\scriptstyle \R\big\{v^\h d(x,h_1)\big\}< \frac{\sigma_v^2 b(x,h_1)}{2} \& \R\big\{v^\h d(x,h_2)\big\}< \frac{\sigma_v^2b(x,h_2)}{2}} \Big\}.
\end{align}
\end{subequations}
Substituting~\eqref{eqn:defintionbxh} and the identity $p(v)=\CN(v;0,\sigma_v^2 I_{n_y})$ into~\eqref{eqn:muxyderivation2}, we get
\begin{align}
I_1=&P\big(v\in\V_1\big|v\sim\CN(v;0,\sigma_v^2)\big)\nonumber\\
&+\exp\bigg[\frac{s_2^2}{\sigma_v^2}\|d(x,h_2)\|^2-s_2b(x,h_2)\bigg]P\big(v\in\V_2\big|v\sim\CN(v;s_2d(x,h_2),\sigma_v^2)\big)\nonumber\\
&+\exp\bigg[\frac{s_1^2}{\sigma_v^2}\|d(x,h_1)\|^2-s_1b(x,h_1)\bigg]P\big(v\in\V_3\big| v\sim\CN(v;s_1d(x,h_1),\sigma_v^2)\big)\nonumber\\
&+\exp\bigg[\frac{1}{\sigma_v^2}\|s_1d(x,h_1)+s_2d(x,h_2)\|^2-s_1b(x,h_1)-s_2b(x,h_2)\bigg]\nonumber\\
&\times P\big(v\in\V_4\big| v\sim\CN(v;s_1d(x,h_1)+s_2d(x,h_2),\sigma_v^2)\big).\label{eqn:muxyderivation3}
\end{align}
We now define the real scalars  $a_1$ and $a_2$ as $a_1\triangleq \R\left\{v^\h d(x,h_1)\right\}$ and $a_2\triangleq \R\left\{v^\h d(x,h_2)\right\}$. Since each of the probabilities on the right hand side of~\eqref{eqn:muxyderivation3} are conditioned on $v$ being distributed with a circular symmetric complex Gaussian density and since $a_1$, $a_2$ are linearly dependent on $v$, we have the vector $a\triangleq [a_1,a_2]^\t$ distributed with a Gaussian density which gives
\begin{align}
I_1=&P\bigg(\substack{a_1\ge\frac{\sigma_v^2 b(x,h_1)}{2}\\a_2\ge\frac{\sigma_v^2b(x,h_2)}{2}}\bigg|  a\sim \N\left(a;\bar{a}'_1(x,h_1,h_2),\frac{\sigma_v^2}{2}\Gamma(x,h_1,h_2)\right)\bigg)\nonumber\\
&+\exp\bigg[\frac{s_2^2}{\sigma_v^2}\|d(x,h_2)\|^2-s_2b(x,h_2)\bigg]P\bigg(\substack{a_1\ge\frac{\sigma_v^2 b(x,h_1)}{2}\\a_2< \frac{\sigma_v^2b(x,h_2)}{2}}\bigg| a\sim\N\left(a;\bar{a}'_2(x,h_1,h_2),\frac{\sigma_v^2}{2}\Gamma(x,h_1,h_2)\right)\bigg)\nonumber\\
&+\exp\bigg[\frac{s_1^2}{\sigma_v^2}\|d(x,h_1)\|^2-s_1b(x,h_1)\bigg]P\bigg(\substack{a_1< \frac{\sigma_v^2 b(x,h_1)}{2}\\a_2\ge\frac{\sigma_v^2b(x,h_2)}{2}}\bigg| a\sim\N\left(a;\bar{a}'_3(x,h_1,h_2),\frac{\sigma_v^2}{2}\Gamma(x,h_1,h_2)\right)\bigg)\nonumber\\
&+\exp\bigg[\frac{1}{\sigma_v^2}\|s_1d(x,h_1)+s_2d(x,h_2)\|^2-s_1b(x,h_1)-s_2b(x,h_2)\bigg]\nonumber\\
&\times P\bigg(\substack{a_1< \frac{\sigma_v^2 b(x,h_1)}{2}\\a_2< \frac{\sigma_v^2b(x,h_2)}{2}}\bigg| a\sim\N\left(a;\bar{a}'_4(x,h_1,h_2),\frac{\sigma_v^2}{2}\Gamma(x,h_1,h_2)\right)\bigg)\label{eqn:muxyderivation4}
\end{align}
where
\begin{subequations}
\begin{align}
\bar{a}'_1(x,h_1,h_2)\triangleq&\left[0,0\right]^\t{},\\
\bar{a}'_2(x,h_1,h_2)\triangleq&s_2\left[\R\big\{d^\h(x,h_1)d(x,h_2)\big\},\|d(x,h_2)\|^2\right]^\t{},\\
\bar{a}'_3(x,h_1,h_2)\triangleq&s_1\left[\|d(x,h_1)\|^2,\R\big\{d^\h(x,h_1)d(x,h_2)\big\}\right]^\t{},\\
\bar{a}'_4(x,h_1,h_2)\triangleq&\left[s_1\|d(x,h_1)\|^2+s_2\R\big\{d^\h(x,h_1)d(x,h_2)\big\}, s_2\|d(x,h_2)\|^2+s_1\R\big\{d^\h(x,h_1)d(x,h_2)\big\}\right]^\t{},
\end{align}
\end{subequations}
\begin{align}
\Gamma(x,h_1,h_2)\triangleq&{\left[\begin{array}{cc}\|d(x,h_1)\|^2&\R\big\{d^\h(x,h_1)d(x,h_2)\big\}\\\R\big\{d^\h(x,h_1)d(x,h_2)\big\} &\|d(x,h_2)\|^2\end{array}\right]}.
\end{align}
Each of the probabilities on the right hand side of~\eqref{eqn:muxyderivation4} can be written using the cumulative distribution function $\Ncdf_2(\cdot,\cdot,\cdot)$ of a bivariate Gaussian random variable to give
\begin{align}
I_1=&\Ncdf_2\left(\left[-\frac{\sigma_v^2 b(x,h_1)}{2},-\frac{\sigma_v^2b(x,h_2)}{2}\right]^\t;\bar{a}_1(x,h_1,h_2),\frac{\sigma_v^2}{2}\Gamma(x,h_1,h_2)\right)\nonumber\\
&+\exp\bigg[\frac{s_2^2}{\sigma_v^2}\|d(x,h_2)\|^2-s_2b(x,h_2)\bigg]\nonumber\\
&\times\Ncdf_2\left(\left[-\frac{\sigma_v^2 b(x,h_1)}{2},\frac{\sigma_v^2b(x,h_2)}{2}\right]^\t;\bar{a}_2(x,h_1,h_2),\frac{\sigma_v^2}{2}\overline{\Gamma}(x,h_1,h_2)\right)\nonumber\\
&+\exp\bigg[\frac{s_1^2}{\sigma_v^2}\|d(x,h_1)\|^2-s_1b(x,h_1)\bigg]\nonumber\\
&\times\Ncdf_2\left(\left[\frac{\sigma_v^2 b(x,h_1)}{2},-\frac{\sigma_v^2b(x,h_2)}{2}\right]^\t;\bar{a}_3(x,h_1,h_2),\frac{\sigma_v^2}{2}\overline{\Gamma}(x,h_1,h_2)\right)\nonumber\\
&+\exp\bigg[\frac{1}{\sigma_v^2}\|s_1d(x,h_1)+s_2d(x,h_2)\|^2-s_1b(x,h_1)-s_2b(x,h_2)\bigg]\nonumber\\
&\times\Ncdf_2\left(\left[\frac{\sigma_v^2 b(x,h_1)}{2},\frac{\sigma_v^2b(x,h_2)}{2}\right]^\t;\bar{a}_4(x,h_1,h_2),\frac{\sigma_v^2}{2}\Gamma(x,h_1,h_2)\right)\label{eqn:muxyderivation5}
\end{align}
where
\begin{subequations}
\begin{align}
\bar{a}_1(x,h_1,h_2)\triangleq&-\bar{a}'_1(x,h_1,h_2),\\
\bar{a}_2(x,h_1,h_2)\triangleq&{\scriptsize\left[\begin{array}{cc}-1&0\\0&1\end{array}\right]}\bar{a}'_2(x,h_1,h_2),\\ \bar{a}_3(x,h_1,h_2)\triangleq&{\scriptsize\left[\begin{array}{cc}1&0\\0&-1\end{array}\right]}\bar{a}'_3(x,h_1,h_2),\\
\bar{a}_4(x,h_1,h_2)\triangleq&\bar{a}'_4(x,h_1,h_2),
\end{align}
\end{subequations}
\begin{align}
\overline{\Gamma}(x,h_1,h_2)\triangleq{\left[\begin{array}{cc}\|d(x,h_1)\|^2&\hspace{-0.5cm}-\R\big\{d^\h(x,h_1)d(x,h_2)\big\}\\-\R\big\{d^\h(x,h_1)d(x,h_2)\big\} &\hspace{-0.5cm}\|d(x,h_2)\|^2\end{array}\right]}.
\end{align}
Using~\eqref{eqn:defintionbxh} in~\eqref{eqn:muxyderivation5}, we get
\begin{align}
I_1=&\Ncdf_2\left(\left[-\frac{\sigma_v^2 b(x,h_1)}{2},-\frac{\sigma_v^2b(x,h_2)}{2}\right]^\t;\bar{a}_1(x,h_1,h_2),\frac{\sigma_v^2}{2}\Gamma(x,h_1,h_2)\right)\nonumber\\
&+\exp\bigg[\frac{s_2^2-s_2}{\sigma_v^2}\|d(x,h_2)\|^2-\frac{s_2}{\sigma_x^2}x^\t h_2-\frac{s_2}{2\sigma_x^2}\|h_2\|^2\bigg]\nonumber\\
&\times\Ncdf_2\left(\left[ -\frac{\sigma_v^2 b(x,h_1)}{2},\frac{\sigma_v^2b(x,h_2)}{2}\right]^\t;\bar{a}_2(x,h_1,h_2),\frac{\sigma_v^2}{2}\overline{\Gamma}(x,h_1,h_2)\right)\nonumber\\
&+\exp\bigg[\frac{s_1^2-s_1}{\sigma_v^2}\|d(x,h_1)\|^2-\frac{s_1}{\sigma_x^2}x^\t h_1-\frac{s_1}{2\sigma_x^2}\|h_1\|^2\bigg]\nonumber\\
&\times\Ncdf_2\left(\left[\frac{\sigma_v^2 b(x,h_1)}{2},-\frac{\sigma_v^2b(x,h_2)}{2}\right]^\t;\bar{a}_3(x,h_1,h_2),\frac{\sigma_v^2}{2}\overline{\Gamma}(x,h_1,h_2)\right)\nonumber\\
&+\exp\bigg[\frac{s_1^2-s_1}{\sigma_v^2}\|d(x,h_1)\|^2+\frac{s_2^2-s_2}{\sigma_v^2}\|d(x,h_2)\|^2+\frac{2s_1s_2}{\sigma_v^2}\R\left\{d^\h(x,h_1)d(x,h_2)\right\}\nonumber\\
&-\frac{s_1h_1+s_2h_2}{\sigma_x^2}x-\frac{s_1}{2\sigma_x^2}\|h_1\|^2-\frac{s_2}{2\sigma_x^2}\|h_2\|^2\bigg]\nonumber\\
&\times\Ncdf_2\left(\left[\frac{\sigma_v^2 b(x,h_1)}{2},\frac{\sigma_v^2b(x,h_2)}{2}\right]^\t; \bar{a}_4(x,h_1,h_2),\frac{\sigma_v^2}{2}\Gamma(x,h_1,h_2)\right)\\
=&\Ncdf_2\left(\left[ -\frac{\sigma_v^2 b(x,h_1)}{2},-\frac{\sigma_v^2b(x,h_2)}{2}\right]^\t;\bar{a}_1(x,h_1,h_2),\frac{\sigma_v^2}{2}\Gamma(x,h_1,h_2)\right)\nonumber\\
&+\exp\bigg[\frac{s_2^2-s_2}{\sigma_v^2}\|d(x,h_2)\|^2+\frac{s_2^2-s_2}{2\sigma_x^2}\|h_2\|^2\bigg]\exp\bigg[-\frac{s_2}{\sigma_x^2}x^\t h_2-\frac{s_2^2}{2\sigma_x^2}\|h_2\|^2\bigg]\nonumber\\
&\times\Ncdf_2\left(\left[ -\frac{\sigma_v^2 b(x,h_1)}{2},\frac{\sigma_v^2b(x,h_2)}{2}\right]^\t;\bar{a}_2(x,h_1,h_2),\frac{\sigma_v^2}{2}\overline{\Gamma}(x,h_1,h_2)\right)\nonumber\\
&+\exp\bigg[\frac{s_1^2-s_1}{\sigma_v^2}\|d(x,h_1)\|^2+\frac{s_1^2-s_1}{2\sigma_x^2}\|h_1\|^2\bigg]\exp\bigg[-\frac{s_1}{\sigma_x^2}x^\t h_1-\frac{s_1^2}{2\sigma_x^2}\|h_1\|^2\bigg]\nonumber\\
&\times\Ncdf_2\left(\left[\frac{\sigma_v^2 b(x,h_1)}{2},-\frac{\sigma_v^2b(x,h_2)}{2}\right]^\t;\bar{a}_3(x,h_1,h_2),\frac{\sigma_v^2}{2}\overline{\Gamma}(x,h_1,h_2)\right)\nonumber\\
&+\exp\bigg[\frac{s_1^2-s_1}{\sigma_v^2}\|d(x,h_1)\|^2+\frac{s_2^2-s_2}{\sigma_v^2}\|d(x,h_2)\|^2+\frac{2s_1s_2}{\sigma_v^2}\R\left\{d^\h(x,h_1)d(x,h_2)\right\}\bigg]\nonumber\\
&\times\exp\bigg[\frac{s_1^2-s_1}{2\sigma_x^2}\|h_1\|^2+\frac{s_2^2-s_2}{2\sigma_x^2}\|h_2\|^2+\frac{s_1s_2}{\sigma_x^2}h_1^\t h_2 \bigg]\exp\bigg[-\frac{s_1h_1+s_2h_2}{\sigma_x^2}x-\frac{1}{2\sigma_x^2}\|s_1h_1+s_2h_2\|^2\bigg]\nonumber\\
&\times\Ncdf_2\left(\left[\frac{\sigma_v^2 b(x,h_1)}{2},\frac{\sigma_v^2b(x,h_2)}{2}\right]^\t; \bar{a}_4(x,h_1,h_2),\frac{\sigma_v^2}{2}\Gamma(x,h_1,h_2)\right).\label{eqn:muxyderivation6}
\end{align}
Substituting the result~\eqref{eqn:muxyderivation6} into~\eqref{eqn:muxyderivation1} and carrying out straightforward algebra, we obtain
\begin{align}
\mu_{y,x}(s_1,s_2,h_1,h_2)=&E_{\N(x;0,\sigma_x^2)}\Bigg[\Ncdf_2\left(\left[\begin{array}{c}-\frac{\sigma_v^2 b(x,h_1)}{2}\\-\frac{\sigma_v^2b(x,h_2)}{2}\end{array}\right];\bar{a}_1(x,h_1,h_2),\frac{\sigma_v^2}{2}\Gamma(x,h_1,h_2)\right)\Bigg]\nonumber\\
&+E_{\N(x;-s_2h_2,\sigma_x^2)}\Bigg[\exp\bigg[\frac{s_2^2-s_2}{\sigma_v^2}\|d(x,h_2)\|^2+\frac{s_2^2-s_2}{2\sigma_x^2}\|h_2\|^2\bigg]\nonumber\\
&\times\Ncdf_2\left(\left[\begin{array}{c}-\frac{\sigma_v^2 b(x,h_1)}{2}\\\frac{\sigma_v^2b(x,h_2)}{2}\end{array}\right];\bar{a}_2(x,h_1,h_2),\frac{\sigma_v^2}{2}\overline{\Gamma}(x,h_1,h_2)\right)\Bigg]\nonumber\\
&+E_{\N(x;-s_1h_1,\sigma_x^2)}\Bigg[\exp\bigg[\frac{s_1^2-s_1}{\sigma_v^2}\|d(x,h_1)\|^2+\frac{s_1^2-s_1}{2\sigma_x^2}\|h_1\|^2\bigg]\nonumber\\
&\times \Ncdf_2\left(\left[\begin{array}{c}\frac{\sigma_v^2 b(x,h_1)}{2}\\-\frac{\sigma_v^2b(x,h_2)}{2}\end{array}\right];\bar{a}_3(x,h_1,h_2),\frac{\sigma_v^2}{2}\overline{\Gamma}(x,h_1,h_2)\right)\Bigg]\nonumber\\
&+E_{\N(x;-(s_1h_1+s_2h_2),\sigma_x^2)}\Bigg[\exp\bigg[\frac{s_1^2-s_1}{\sigma_v^2}\|d(x,h_1)\|^2+\frac{s_2^2-s_2}{\sigma_v^2}\|d(x,h_2)\|^2\nonumber\\
&+\frac{2s_1s_2}{\sigma_v^2}\R\left\{d^\h(x,h_1)d(x,h_2)\right\}\bigg]\exp\bigg[\frac{s_1^2-s_1}{2\sigma_x^2}\|h_1\|^2+\frac{s_2^2-s_2}{2\sigma_x^2}\|h_2\|^2+\frac{s_1s_2}{\sigma_x^2}h_1^\t h_2 \bigg]\nonumber\\
&\times \Ncdf_2\left(\left[\begin{array}{c}\frac{\sigma_v^2 b(x,h_1)}{2}\\\frac{\sigma_v^2b(x,h_2)}{2}\end{array}\right];\bar{a}_4(x,h_1,h_2),\frac{\sigma_v^2}{2}\Gamma(x,h_1,h_2)\right)\Bigg].\label{eqn:formulamuxy}
\end{align}

\section{Calculation of $\mu_{y|x}(s_1,s_2,h_1,h_2,x)$ and $\mu_{x}(s_1,s_2,h_1,h_2)$}
We first calculate $L_1(\cdot,\cdot,\cdot)$ and $L_2(\cdot,\cdot)$ as follows.
\begin{align}
L_1(y,x+h,x)=&\frac{\exp\left[-\frac{1}{\sigma_v^2}\|y-g(x+h)\|^2\right]}{\exp\left[-\frac{1}{\sigma_v^2}\|y-g(x)\|^2\right]}\\
=&\exp\left[\frac{2}{\sigma_v^2}\R\left\{y^\h d(x,h)\right\}-\frac{1}{\sigma_v^2}\|g(x+h)\|^2+\frac{1}{\sigma_v^2}\|g(x)\|^2\right]\\
=&\exp\left[\frac{2}{\sigma_v^2}\R\left\{v^\h d(x,h)\right\}-\frac{1}{\sigma_v^2}\|d(x,h)\|^2\right]\\
=&\exp\left[\frac{2}{\sigma_v^2}\R\left\{v^\h d(x,h)\right\}-b_1(x,h)\right],\\
L_2(x+h,x)=&\frac{\exp\left[-\frac{1}{2\sigma_x^2}\|x+h\|^2\right]}{\exp\left[-\frac{1}{2\sigma_x^2}\|x\|^2\right]}\\
=&\exp\left[-\frac{1}{\sigma_x^2}x^\t h-\frac{1}{2\sigma_x^2}\|h\|^2\right]=\exp\left[-b_2(x,h)\right],
\end{align}
where
\begin{align}
b_1(x,h)\triangleq&\frac{1}{\sigma_v^2}\|d(x,h)\|^2,\label{eqn:defintionb1xh}\\
b_2(x,h)\triangleq&\frac{1}{\sigma_x^2}x^\t h+\frac{1}{2\sigma_x^2}\|h\|^2.\label{eqn:defintionb2xh}
\end{align}
Then the function $\mu_{y|x}(\cdot,\cdot,\cdot,\cdot,\cdot)$ is given as
\begin{align}
\mu_{y|x}&(s_1,s_2,h_1,h_2,x)\nonumber\\
=&\int \min\left(\exp\left[\frac{2s_1}{\sigma_v^2}\R\left\{v^\h d(x,h_1)\right\}-s_1b_1(x,h_1)\right],1\right)\nonumber\\
&\times\min\left(\exp\left[\frac{2s_2}{\sigma_v^2}\R\left\{v^\h d(x,h_2)\right\}-s_2b_1(x,h_2)\right],1\right)p(v)\d v\\
=&\int_{\V_1}p(v)\d v\nonumber\\
&+\int_{\V_2}\exp\left[\frac{2s_2}{\sigma_v^2}\R\left\{v^\h d(x,h_2)\right\}-s_2b_1(x,h_2)\right]p(v)\d v\nonumber\\
&+\int_{\V_3} \exp\left[\frac{2s_1}{\sigma_v^2}\R\left\{v^\h d(x,h_1)\right\}-s_1b_1(x,h_1)\right]p(v)\d v\nonumber\\
&+\int_{\V_4} \exp\left[\frac{2s_1}{\sigma_v^2}\R\left\{v^\h d(x,h_1)\right\}-s_1b_1(x,h_1)\right]\exp\left[\frac{2s_2}{\sigma_v^2}\R\left\{v^\h d(x,h_2)\right\}-s_2b_1(x,h_2)\right]p(v)\d v\label{eqn:muxgyderivation2}
\end{align}
where the sets $\V_1$, $\V_2$, $\V_3$ and $\V_4$ are defined as follows.
\begin{subequations}
\begin{align}
\V_1&\triangleq\Big\{v\in\mathbb{C}^{n_y}\Big|{\scriptstyle \R\big\{v^\h d(x,h_1)\big\}\ge\frac{\sigma_v^2 b_1(x,h_1)}{2} \& \R\big\{v^\h d(x,h_2)\big\}\ge\frac{\sigma_v^2b_1(x,h_2)}{2}} \Big\},\\
\V_2&\triangleq\Big\{v\in\mathbb{C}^{n_y}\Big|{\scriptstyle \R\big\{v^\h d(x,h_1)\big\}\ge\frac{\sigma_v^2 b_1(x,h_1)}{2} \& \R\big\{v^\h d(x,h_2)\big\}< \frac{\sigma_v^2b_1(x,h_2)}{2}} \Big\},\\
\V_3&\triangleq\Big\{v\in\mathbb{C}^{n_y}\Big|{\scriptstyle \R\big\{v^\h d(x,h_1)\big\}< \frac{\sigma_v^2 b_1(x,h_1)}{2} \& \R\big\{v^\h d(x,h_2)\big\}\ge\frac{\sigma_v^2b_1(x,h_2)}{2}} \Big\},\\
\V_4&\triangleq\Big\{v\in\mathbb{C}^{n_y}\Big|{\scriptstyle \R\big\{v^\h d(x,h_1)\big\}< \frac{\sigma_v^2 b_1(x,h_1)}{2} \& \R\big\{v^\h d(x,h_2)\big\}< \frac{\sigma_v^2b_1(x,h_2)}{2}} \Big\}.
\end{align}
\end{subequations}
Substituting~\eqref{eqn:defintionb1xh} and the identity $p(v)=\CN(v;0,\sigma_v^2 I_{n_y})$ into~\eqref{eqn:muxgyderivation2}, we get

\begin{align}
\mu_{y|x}&(s_1,s_2,h_1,h_2,x)\nonumber\\
=&P\Bigg(\substack{\R\left\{v^\h d(x,h_1)\right\}\ge\frac{\sigma_v^2 b_1(x,h_1)}{2}\\\R\left\{v^\h d(x,h_2)\right\}\ge\frac{\sigma_v^2 b_1(x,h_2)}{2}}\Bigg|v\sim\CN(v;0,\sigma_v^2I_{n_y})\Bigg)\\
+&\exp\left[\frac{s_2^2-s_2}{\sigma_v^2}\|d(x,h_2)\|^2\right]P\Bigg(\substack{\R\left\{v^\h d(x,h_1)\right\}\ge\frac{\sigma_v^2 b_1(x,h_1)}{2}\\\R\left\{v^\h d(x,h_2)\right\}<\frac{\sigma_v^2 b_1(x,h_2)}{2}}\Bigg|v\sim\CN(v;s_2d(x,h_2),\sigma_v^2I_{n_y})\Bigg)\nonumber\\
+&\exp\left[\frac{s_1^2-s_1}{\sigma_v^2}\|d(x,h_1)\|^2\right]P\Bigg(\substack{\R\left\{v^\h d(x,h_1)\right\}<\frac{\sigma_v^2 b_1(x,h_1)}{2}\\\R\left\{v^\h d(x,h_2)\right\}\ge\frac{\sigma_v^2 b_1(x,h_2)}{2}}\Bigg|v\sim\CN(v;s_1d(x,h_1),\sigma_v^2I_{n_y})\Bigg)\nonumber\\
&+\exp\left[\frac{s_1^2-s_1}{\sigma_v^2}\|d(x,h_1)\|^2+\frac{s_2^2-s_2}{\sigma_v^2}\|d(x,h_2)\|^2+\frac{2s_1s_2}{\sigma_v^2}\R\{d^\h(x,h_1)d(x,h_2)\}\right]\nonumber\\
&\times P\Bigg(\substack{\R\left\{v^\h d(x,h_1)\right\}<\frac{\sigma_v^2 b_1(x,h_1)}{2}\\\R\left\{v^\h d(x,h_2)\right\}<\frac{\sigma_v^2 b_1(x,h_2)}{2}}\Bigg|v\sim\CN(v;s_1d(x,h_1)+s_2d(x,h_2),\sigma_v^2I_{n_y})\Bigg).
\end{align}
Continuing in the same way as in Section~\ref{sec:muyx}, we get
\begin{align}
\mu_{y|x}(s_1,s_2,h_1,h_2,x)=&\Ncdf_2\left(\left[\begin{array}{c}-\frac{\sigma_v^2 b_1(x,h_1)}{2}\\-\frac{\sigma_v^2b_1(x,h_2)}{2}\end{array}\right];\bar{a}_1(x,h_1,h_2),\frac{\sigma_v^2}{2}\Gamma(x,h_1,h_2)\right)\nonumber\\
&+\exp\bigg[\frac{s_2^2-s_2}{\sigma_v^2}\|d(x,h_2)\|^2\bigg]\nonumber\\
&\times\Ncdf_2\left(\left[\begin{array}{c}-\frac{\sigma_v^2 b_1(x,h_1)}{2}\\\frac{\sigma_v^2b_1(x,h_2)}{2}\end{array}\right];\bar{a}_2(x,h_1,h_2),\frac{\sigma_v^2}{2}\overline{\Gamma}(x,h_1,h_2)\right)\nonumber\\
&+\exp\bigg[\frac{s_1^2-s_1}{\sigma_v^2}\|d(x,h_1)\|^2\nonumber\\
&\times\Ncdf_2\left(\left[\begin{array}{c}\frac{\sigma_v^2 b_1(x,h_1)}{2}\\-\frac{\sigma_v^2b_1(x,h_2)}{2}\end{array}\right];\bar{a}_3(x,h_1,h_2),\frac{\sigma_v^2}{2}\overline{\Gamma}(x,h_1,h_2)\right)\bigg]\nonumber\\
&+\exp\bigg[\frac{s_1^2-s_1}{\sigma_v^2}\|d(x,h_1)\|^2+\frac{s_2^2-s_2}{\sigma_v^2}\|d(x,h_2)\|^2+\frac{2s_1s_2}{\sigma_v^2}\R\left\{d^\h(x,h_1)d(x,h_2)\right\}\bigg]\nonumber\\
&\times\Ncdf_2\left(\left[\begin{array}{c}\frac{\sigma_v^2 b_1(x,h_1)}{2}\\\frac{\sigma_v^2b_1(x,h_2)}{2}\end{array}\right];\bar{a}_4(x,h_1,h_2),\frac{\sigma_v^2}{2}\Gamma(x,h_1,h_2)\right).\label{eqn:formulamuxgy}
\end{align}
Similarly the function $\mu_{x}(\cdot,\cdot,\cdot,\cdot)$ is given as
\begin{align}
\mu_{x}&(s_1,s_2,h_1,h_2)\nonumber\\
=&\int \min\left(\exp\left[-s_1b_2(x,h_1)\right],1\right)\min\left(\exp\left[-s_2b_2(x,h_2)\right],1\right)p(x)\d x\\
=&\int_{\substack{b_2(x,h_1)<0\\b_2(x,h_2)<0}}p(x)\d x+\int_{\substack{b_2(x,h_1)<0\\b_2(x,h_2)\ge0}} \exp\left[-s_2b_2(x,h_2)\right]p(x)\d x\nonumber\\
&+\int_{\substack{b_2(x,h_1)\ge0\\b_2(x,h_2)<0}} \exp\left[-s_1b_2(x,h_1)\right]p(x)\d x+\int_{\substack{b_2(x,h_1)\ge0\\b_2(x,h_2)\ge0}} \exp\left[-s_1b_2(x,h_1)\right]\exp\left[-s_2b_2(x,h_2)\right]p(x)\d x\\
=&P\Big(\substack{b_2(x,h_1)<0\\b_2(x,h_2)<0}\Big| x\sim\N(x;0;\sigma_x^2)\Big)\nonumber\\
&+\exp\left[\frac{s_2^2-s_2}{2\sigma_x^2}\|h_2\|^2\right]P\Big(\substack{b_2(x,h_1)<0\\b_2(x,h_2)\ge0}\Big|x\sim\N(x;-s_2h_2,\sigma_x^2)\Big)\nonumber\\
&+\exp\left[\frac{s_1^2-s_1}{2\sigma_x^2}\|h_1\|^2\right]P\Big(\substack{b_2(x,h_1)\ge0\\b_2(x,h_2)<0}\Big|x\sim\N(x;-s_1h_1,\sigma_x^2)\Big)\nonumber\\
&+\exp\left[\frac{s_1^2-s_1}{2\sigma_x^2}\|h_1\|^2+\frac{s_2^2-s_2}{2\sigma_x^2}\|h_2\|^2+\frac{s_1s_2}{\sigma_x^2}h_1^\t h_2\right] P\Big(\substack{b_2(x,h_1)\ge0\\b_2(x,h_2)\ge0}\Big|x\sim\N(x;-s_1h_1-s_2h_2,\sigma_x^2)\Big)\\
=&P\Bigg(\begin{array}{c}a_1<0\\a_2<0\end{array}\Bigg| a\sim\N\left(a;\tilde{a}'_1(h_1,h_2),\frac{1}{\sigma_x^2}\Lambda(h_1,h_2)\right)\Bigg)\nonumber\\
&+\exp\left[\frac{s_2^2-s_2}{2\sigma_x^2}\|h_2\|^2\right]P\Bigg(\begin{array}{c}a_1<0\\a_2\ge0\end{array}\Bigg| a\sim\N\left(a;\tilde{a}'_2(h_1,h_2),\frac{1}{\sigma_x^2}\Lambda(h_1,h_2)\right)\Bigg)\nonumber\\
&+\exp\left[\frac{s_1^2-s_1}{2\sigma_x^2}\|h_1\|^2\right]P\Bigg(\begin{array}{c}a_1\ge0\\a_2<0\end{array}\Bigg| a\sim\N\left(a;\tilde{a}'_3(h_1,h_2),\frac{1}{\sigma_x^2}\Lambda(h_1,h_2)\right)\Bigg)\nonumber\\
&+\exp\left[\frac{s_1^2-s_1}{2\sigma_x^2}\|h_1\|^2+\frac{s_2^2-s_2}{2\sigma_x^2}\|h_2\|^2+\frac{s_1s_2}{\sigma_x^2}h_1^\t h_2\right]\nonumber\\
&\times P\Bigg(\begin{array}{c}a_1\ge0\\a_2\ge0\end{array}\Bigg| a\sim\N\left(a;\tilde{a}'_4(h_1,h_2),\frac{1}{\sigma_x^2}\Lambda(h_1,h_2)\right)\Bigg)\label{eqn:muxderivation1}
\end{align}
where
\begin{subequations}
\begin{align}
\tilde{a}'_1(h_1,h_2)\triangleq&\frac{1}{2\sigma_x^2}\left[\|h_1\|^2,\|h_2\|^2\right]^\t{},\\
\tilde{a}'_2(h_1,h_2)\triangleq&\frac{1}{2\sigma_x^2}\left[-2s_2h_1^\t h_2+\|h_1\|^2,(1-2s_2)\|h_2\|^2\right]^\t{},\\
\tilde{a}'_3(h_1,h_2)\triangleq&\frac{1}{2\sigma_x^2}\left[(1-2s_1)\|h_1\|^2,-2s_1h_1^\t h_2+\|h_2\|^2\right]^\t{},\\
\tilde{a}'_4(h_1,h_2)\triangleq&\frac{1}{2\sigma_x^2}\left[\begin{array}{c}(1-2s_1)\|h_1\|^2-2s_2h_1^\t h_2 \\-2s_1h_1^\t h_2+(1-2s_2)\|h_2\|^2\end{array}\right],
\end{align}
\end{subequations}
\begin{align}
\Lambda(h_1,h_2)\triangleq&\left[\begin{array}{cc}\|h_1\|^2& h_1^\t h_2\\h_1^\t h_2 &\|h_2\|^2 \end{array}\right].
\end{align}
Each of the probabilities on the right hand side of~\eqref{eqn:muxderivation1} can be written using the cumulative distribution function $\Ncdf_2(\cdot,\cdot,\cdot)$ as
\begin{align}
\mu_{x}&(s_1,s_2,h_1,h_2)\nonumber\\
=&\Ncdf_2\Bigg(\left[\begin{array}{c}0\\0\end{array}\right];\tilde{a}_1(h_1,h_2),\frac{1}{\sigma_x^2}\Lambda(h_1,h_2)\Bigg)\nonumber\\
&+\exp\left[\frac{s_2^2-s_2}{2\sigma_x^2}\|h_2\|^2\right]\Ncdf_2\Bigg(\left[\begin{array}{c}0\\0\end{array}\right];\tilde{a}_2(h_1,h_2),\frac{1}{\sigma_x^2}\overline{\Lambda}(h_1,h_2)\Bigg)\nonumber\\
&+\exp\left[\frac{s_1^2-s_1}{2\sigma_x^2}\|h_1\|^2\right]\Ncdf_2\Bigg(\left[\begin{array}{c}0\\0\end{array}\right];\tilde{a}_3(h_1,h_2),\frac{1}{\sigma_x^2}\overline{\Lambda}(h_1,h_2)\Bigg)\nonumber\\
&+\exp\left[\frac{s_1^2-s_1}{2\sigma_x^2}\|h_1\|^2+\frac{s_2^2-s_2}{2\sigma_x^2}\|h_2\|^2+\frac{s_1s_2}{\sigma_x^2}h_1^\t h_2\right]\Ncdf_2\Bigg(\left[\begin{array}{c}0\\0\end{array}\right];\tilde{a}_4(h_1,h_2),\frac{1}{\sigma_x^2}\Lambda(h_1,h_2)\Bigg)\label{eqn:formulamux}
\end{align}
where
\begin{subequations}
\begin{align}
\tilde{a}_1(h_1,h_2)\triangleq&\tilde{a}'_1(h_1,h_2),\\
\tilde{a}_2(h_1,h_2)\triangleq&{\scriptsize\left[\begin{array}{cc}1&0\\0&-1\end{array}\right]}\tilde{a}'_2(h_1,h_2),\\ \tilde{a}_3(h_1,h_2)\triangleq&{\scriptsize\left[\begin{array}{cc}-1&0\\0&1\end{array}\right]}\tilde{a}'_3(h_1,h_2),\\
\tilde{a}_4(h_1,h_2)\triangleq&-\tilde{a}'_4(h_1,h_2),
\end{align}
\end{subequations}
\begin{align}
\overline{\Lambda}(h_1,h_2)\triangleq&\left[\begin{array}{cc}\|h_1\|^2& -h_1^\t h_2\\-h_1^\t h_2 &\|h_2\|^2 \end{array}\right].
\end{align}
\section{Special Cases} In this section, we are going to find the expressions for the following quantities one by one: $\mu_{y,x}(1,1,h,h)$, $\mu_{y,x}(1,0,h,h)$,
$\mu_{y|x}(1,1,h,h,x)$, $\mu_{y|x}(1,0,h,h,x)$, $\mu_{x}(1,1,h,h)$, $\mu_{x}(1,0,h,h)$.
\begin{itemize}
\item \textbf{$\mu_{y,x}(1,1,h,h)$}: Substituting $s_1=s_2=1$ and $h_1=h_2=h$ in~\eqref{eqn:formulamuxy}, we get
\begin{align}
\mu_{y,x}(1,1,h,h)=&E_{\N(x;0,\sigma_x^2)}\Bigg[\Ncdf_2\left(\left[\begin{array}{c}-\frac{\sigma_v^2 b(x,h)}{2}\\-\frac{\sigma_v^2b(x,h)}{2}\end{array}\right];\bar{a}_1^{1,1}(x,h,h),\frac{\sigma_v^2}{2}\Gamma(x,h,h)\right)\Bigg]\nonumber\\
&+E_{\N(x;-h,\sigma_x^2)}\Bigg[\Ncdf_2\left(\left[\begin{array}{c}-\frac{\sigma_v^2 b(x,h)}{2}\\\frac{\sigma_v^2b(x,h)}{2}\end{array}\right];\bar{a}_2^{1,1}(x,h,h),\frac{\sigma_v^2}{2}\overline{\Gamma}(x,h,h)\right)\Bigg]\nonumber\\
&+E_{\N(x;-h,\sigma_x^2)}\Bigg[\Ncdf_2\left(\left[\begin{array}{c}\frac{\sigma_v^2 b(x,h)}{2}\\-\frac{\sigma_v^2b(x,h)}{2}\end{array}\right];\bar{a}_3^{1,1}(x,h,h),\frac{\sigma_v^2}{2}\overline{\Gamma}(x,h,h)\right)\Bigg]\nonumber\\
&+E_{\N(x;-2h,\sigma_x^2)}\Bigg[\exp\bigg[\frac{2}{\sigma_v^2}\|d(x,h)\|^2\bigg]\exp\bigg[\frac{1}{\sigma_x^2}\|h\|^2 \bigg]\nonumber\\
&\times \Ncdf_2\left(\left[\begin{array}{c}\frac{\sigma_v^2 b(x,h)}{2}\\\frac{\sigma_v^2b(x,h)}{2}\end{array}\right];\bar{a}_4^{1,1}(x,h,h),\frac{\sigma_v^2}{2}\Gamma(x,h,h)\right)\Bigg]\label{eqn:formulamuxyspecial1}
\end{align}
where
\begin{subequations}
\begin{align}
\bar{a}_1^{1,1}(x,h,h)\triangleq&\left[0,0\right]^\t,\\
\bar{a}_2^{1,1}(x,h,h)\triangleq&\left[-\|d(x,h)\|^2,\|d(x,h)\|^2\right]^\t,\\
\bar{a}_3^{1,1}(x,h,h)\triangleq&\left[\|d(x,h)\|^2,-\|d(x,h)\|^2\right]^\t,\\
\bar{a}_4^{1,1}(x,h,h)\triangleq&2\left[\|d(x,h)\|^2, \|d(x,h)\|^2\right]^\t,
\end{align}
\end{subequations}
and
\begin{subequations}
\begin{align}
\Gamma(x,h,h)\triangleq&{\left[\begin{array}{cc}\|d(x,h)\|^2&\|d(x,h)\|^2\\\|d(x,h)\|^2 &\|d(x,h)\|^2\end{array}\right]},\\
\overline{\Gamma}(x,h,h)\triangleq&{\left[\begin{array}{cc}\|d(x,h)\|^2&-\|d(x,h)\|^2\\-\|d(x,h)\|^2 &\hspace{-0.5cm}\|d(x,h)\|^2\end{array}\right]}.
\end{align}
\end{subequations}
Noting that $\Gamma(x,h,h)$ and $\overline{\Gamma}(x,h,h)$ are singular, the quantities related to the bivariate cumulative distribution function $\Ncdf_2(\cdot)$ can be written as
\begin{subequations}
\begin{align}
\Ncdf_2\left(\left[\begin{array}{c}-\frac{\sigma_v^2 b(x,h)}{2}\\-\frac{\sigma_v^2b(x,h)}{2}\end{array}\right];\bar{a}_1(x,h,h),\frac{\sigma_v^2}{2}\Gamma(x,h,h)\right)=&
\Ncdf_1\left(-\frac{\sigma_v^2 b(x,h)}{2};0,\frac{\sigma_v^2}{2}\|d(x,h)\|^2\right),\\
\Ncdf_2\left(\left[\begin{array}{c}-\frac{\sigma_v^2 b(x,h)}{2}\\\frac{\sigma_v^2b(x,h)}{2}\end{array}\right];\bar{a}_2(x,h,h),\frac{\sigma_v^2}{2}\overline{\Gamma}(x,h,h)\right)=&0,\\
\Ncdf_2\left(\left[\begin{array}{c}\frac{\sigma_v^2 b(x,h)}{2}\\-\frac{\sigma_v^2b(x,h)}{2}\end{array}\right];\bar{a}_3(x,h,h),\frac{\sigma_v^2}{2}\overline{\Gamma}(x,h,h)\right)=&0,\\
\Ncdf_2\left(\left[\begin{array}{c}\frac{\sigma_v^2 b(x,h)}{2}\\\frac{\sigma_v^2b(x,h)}{2}\end{array}\right];\bar{a}_4(x,h,h),\frac{\sigma_v^2}{2}\Gamma(x,h,h)\right)=&\Ncdf_1\left(\frac{\sigma_v^2 b(x,h)}{2};2\|d(x,h)\|^2,\frac{\sigma_v^2}{2}\|d(x,h)\|^2\right),
\end{align}
\end{subequations}
where $\Ncdf_1(\cdot)$ is the cumulative distribution function for the univariate Gaussian density, which gives
\begin{align}
\mu_{y,x}(1,1,h,h)=&E_{\N(x;0,\sigma_x^2)}\Bigg[\Ncdf_1\left(-\frac{\sigma_v^2 b(x,h)}{2};0,\frac{\sigma_v^2}{2}\|d(x,h)\|^2\right)\Bigg]\nonumber\\
&+E_{\N(x;-2h,\sigma_x^2)}\Bigg[\exp\bigg[\frac{2}{\sigma_v^2}\|d(x,h)\|^2\bigg]\exp\bigg[\frac{1}{\sigma_x^2}\|h\|^2 \bigg]\nonumber\\
&\times\Ncdf_1\left(\frac{\sigma_v^2 b(x,h)}{2};2\|d(x,h)\|^2,\frac{\sigma_v^2}{2}\|d(x,h)\|^2\right)\Bigg].\label{eqn:formulamuxyspecial2}
\end{align}
\item \textbf{$\mu_{y,x}(1,0,h,h)$}: Substituting $s_1=1$, $s_2=0$ and $h_1=h_2=h$ in~\eqref{eqn:formulamuxy}, we get
\begin{align}
\mu_{y,x}(1,0,h,h)=&E_{\N(x;0,\sigma_x^2)}\Bigg[\Ncdf_2\left(\left[\begin{array}{c}-\frac{\sigma_v^2 b(x,h)}{2}\\-\frac{\sigma_v^2b(x,h)}{2}\end{array}\right];\bar{a}_1^{1,0}(x,h,h),\frac{\sigma_v^2}{2}\Gamma(x,h,h)\right)\Bigg]\nonumber\\
&+E_{\N(x;0,\sigma_x^2)}\Bigg[\Ncdf_2\left(\left[\begin{array}{c}-\frac{\sigma_v^2 b(x,h)}{2}\\\frac{\sigma_v^2b(x,h)}{2}\end{array}\right];\bar{a}_2^{1,0}(x,h,h),\frac{\sigma_v^2}{2}\overline{\Gamma}(x,h,h)\right)\Bigg]\nonumber\\
&+E_{\N(x;-h,\sigma_x^2)}\Bigg[\Ncdf_2\left(\left[\begin{array}{c}\frac{\sigma_v^2 b(x,h)}{2}\\-\frac{\sigma_v^2b(x,h)}{2}\end{array}\right];\bar{a}_3^{1,0}(x,h,h),\frac{\sigma_v^2}{2}\overline{\Gamma}(x,h,h)\right)\Bigg]\nonumber\\
&+E_{\N(x;-h,\sigma_x^2)}\Bigg[\Ncdf_2\left(\left[\begin{array}{c}\frac{\sigma_v^2 b(x,h)}{2}\\\frac{\sigma_v^2b(x,h)}{2}\end{array}\right];\bar{a}_4^{1,0}(x,h,h),\frac{\sigma_v^2}{2}\Gamma(x,h,h)\right)\Bigg]\label{eqn:formulamuxyspecial3}
\end{align}
where
\begin{subequations}
\begin{align}
\bar{a}_1^{1,0}(x,h,h)\triangleq&\left[0,0\right]^\t{},\\
\bar{a}_2^{1,0}(x,h,h)\triangleq&\left[0,0\right]^\t{},\\
\bar{a}_3^{1,0}(x,h,h)\triangleq&\left[\|d(x,h)\|^2,-\|d(x,h)\|^2\right]^\t{},\\
\bar{a}_4^{1,0}(x,h,h)\triangleq&\left[\|d(x,h)\|^2,\|d(x,h)\|^2\right]^\t{}.
\end{align}
\end{subequations}
Again due to the singularity of $\Gamma(x,h,h)$ and $\overline{\Gamma}(x,h,h)$, we have
\begin{align}
\mu_{y,x}(1,0,h,h)=&E_{\N(x;0,\sigma_x^2)}\Bigg[\Ncdf_1\left(-\frac{\sigma_v^2 b(x,h)}{2};0,\frac{\sigma_v^2}{2}\|d(x,h)\|^2\right)\Bigg]\nonumber\\
&+E_{\N(x;-h,\sigma_x^2)}\Bigg[\Ncdf_1\left(\frac{\sigma_v^2 b(x,h)}{2};\|d(x,h)\|^2,\frac{\sigma_v^2}{2}\|d(x,h)\|^2\right)\Bigg].\label{eqn:formulamuxyspecial4}
\end{align}
\item \textbf{$\mu_{y|x}(1,1,h,h,x)$}: Substituting $s_1=s_2=1$ and $h_1=h_2=h$ in~\eqref{eqn:formulamuxgy}, we get
\begin{align}
\mu_{y|x}(1,1,h,h,x)=&\Ncdf_2\left(\left[\begin{array}{c}-\frac{\sigma_v^2 b_1(x,h)}{2}\\-\frac{\sigma_v^2b_1(x,h)}{2}\end{array}\right];\bar{a}_1^{1,1}(x,h,h),\frac{\sigma_v^2}{2}\Gamma(x,h,h)\right)\nonumber\\
&+\Ncdf_2\left(\left[\begin{array}{c}-\frac{\sigma_v^2 b_1(x,h)}{2}\\\frac{\sigma_v^2b_1(x,h)}{2}\end{array}\right];\bar{a}_2^{1,1}(x,h,h),\frac{\sigma_v^2}{2}\overline{\Gamma}(x,h,h)\right)\nonumber\\
&+\Ncdf_2\left(\left[\begin{array}{c}\frac{\sigma_v^2 b_1(x,h)}{2}\\-\frac{\sigma_v^2b_1(x,h)}{2}\end{array}\right];\bar{a}_3^{1,1}(x,h,h),\frac{\sigma_v^2}{2}\overline{\Gamma}(x,h,h)\right)\bigg]\nonumber\\
&+\exp\bigg[\frac{2}{\sigma_v^2}\|d(x,h)\|^2\bigg]\Ncdf_2\left(\left[\begin{array}{c}\frac{\sigma_v^2 b_1(x,h)}{2}\\\frac{\sigma_v^2b_1(x,h)}{2}\end{array}\right];\bar{a}_4(x,h,h),\frac{\sigma_v^2}{2}\Gamma(x,h,h)\right).\label{eqn:formulamuxgyspecial1}
\end{align}
Using the singularity of $\Gamma(x,h,h)$ and $\overline{\Gamma}(x,h,h)$, we obtain
\begin{align}
\mu_{y|x}(1,1,h,h,x)=&\Ncdf_1\left(-\frac{\sigma_v^2 b_1(x,h)}{2};0,\frac{\sigma_v^2}{2}\|d(x,h)\|^2\right)\nonumber\\
&+\exp\bigg[\frac{2}{\sigma_v^2}\|d(x,h)\|^2\bigg]\Ncdf_1\left(\frac{\sigma_v^2 b_1(x,h)}{2};2\|d(x,h)\|^2,\frac{\sigma_v^2}{2}\|d(x,h)\|^2\right).\label{eqn:formulamuxgyspecial2}
\end{align}
\item \textbf{$\mu_{y|x}(1,0,h,h,x)$}: Substituting $s_1=1$, $s_2=0$ and $h_1=h_2=h$ in~\eqref{eqn:formulamuxgy}, we get
\begin{align}
\mu_{y|x}(1,0,h,h,x)=&\Ncdf_2\left(\left[\begin{array}{c}-\frac{\sigma_v^2 b_1(x,h)}{2}\\-\frac{\sigma_v^2b_1(x,h)}{2}\end{array}\right];\bar{a}_1^{1,0}(x,h,h),\frac{\sigma_v^2}{2}\Gamma(x,h,h)\right)\nonumber\\
&+\Ncdf_2\left(\left[\begin{array}{c}-\frac{\sigma_v^2 b_1(x,h)}{2}\\\frac{\sigma_v^2b_1(x,h)}{2}\end{array}\right];\bar{a}_2^{1,0}(x,h,h),\frac{\sigma_v^2}{2}\overline{\Gamma}(x,h,h)\right)\nonumber\\
&+\Ncdf_2\left(\left[\begin{array}{c}\frac{\sigma_v^2 b_1(x,h)}{2}\\-\frac{\sigma_v^2b_1(x,h)}{2}\end{array}\right];\bar{a}_3^{1,0}(x,h,h),\frac{\sigma_v^2}{2}\overline{\Gamma}(x,h,h)\right)\bigg]\nonumber\\
&+\Ncdf_2\left(\left[\begin{array}{c}\frac{\sigma_v^2 b_1(x,h)}{2}\\\frac{\sigma_v^2b_1(x,h)}{2}\end{array}\right];\bar{a}_4^{1,0}(x,h,h),\frac{\sigma_v^2}{2}\Gamma(x,h,h)\right).\label{eqn:formulamuxgyspecial3}
\end{align}
Using the singularity of $\Gamma(x,h,h)$ and $\overline{\Gamma}(x,h,h)$, we obtain
\begin{align}
\mu_{y|x}(1,0,h,h,x)=&\Ncdf_1\left(-\frac{\sigma_v^2 b_1(x,h)}{2};0,\frac{\sigma_v^2}{2}\|d(x,h)\|^2\right)\nonumber\\
&+\Ncdf_1\left(\frac{\sigma_v^2 b_1(x,h)}{2};\|d(x,h)\|^2,\frac{\sigma_v^2}{2}\|d(x,h)\|^2\right).\label{eqn:formulamuxgyspecial4}
\end{align}
If we now substitute $b_1(\cdot,\cdot)$ from~\eqref{eqn:defintionb1xh} into~\eqref{eqn:formulamuxgyspecial4}, we get
\begin{align}
\mu_{y|x}(1,0,h,h,x)=&\Ncdf_1\left(-\frac{\|d(x,h)\|^2}{2};0,\frac{\sigma_v^2}{2}\|d(x,h)\|^2\right)\nonumber\\
&+\Ncdf_1\left(\frac{\|d(x,h)\|^2}{2};\|d(x,h)\|^2,\frac{\sigma_v^2}{2}\|d(x,h)\|^2\right)\\
=&2\Ncdf_1\left(-\frac{\|d(x,h)\|^2}{2};0,\frac{\sigma_v^2}{2}\|d(x,h)\|^2\right)\\
=&1-\erf\left(\frac{\|d(x,h)\|}{2\sigma_v}\right)\label{eqn:formulamuxgyspecial5}
\end{align}
where we used the identity
\begin{align}
\Ncdf_1\left(\xi;\bar{\xi},\sigma_\xi^2\right)=\frac{1}{2}\left(1+\erf\left(\frac{\xi-\bar{\xi}}{\sqrt{2}\sigma_\xi}\right)\right).
\end{align}
\item \textbf{$\mu_{x}(1,1,h,h)$}: Substituting $s_1=s_2=1$ and $h_1=h_2=h$ in~\eqref{eqn:formulamux}, we get
\begin{align}
\mu_{x}(1,1,h,h)=&\Ncdf_2\Bigg(\left[\begin{array}{c}0\\0\end{array}\right];\tilde{a}_1^{1,1}(h,h),\frac{1}{\sigma_x^2}\Lambda(h_1,h_2)\Bigg)\nonumber\\
&+\Ncdf_2\Bigg(\left[\begin{array}{c}0\\0\end{array}\right];\tilde{a}_2^{1,1}(h,h),\frac{1}{\sigma_x^2}\overline{\Lambda}(h,h)\Bigg)\nonumber\\
&+\Ncdf_2\Bigg(\left[\begin{array}{c}0\\0\end{array}\right];\tilde{a}_3^{1,1}(h,h),\frac{1}{\sigma_x^2}\overline{\Lambda}(h,h)\Bigg)\nonumber\\
&+\exp\left[\frac{1}{\sigma_x^2}\|h\|^2\right]\Ncdf_2\Bigg(\left[\begin{array}{c}0\\0\end{array}\right];\tilde{a}_4^{1,1}(h,h),\frac{1}{\sigma_x^2}\Lambda(h,h)\Bigg)\label{eqn:formulamuxspecial1}
\end{align}
where
\begin{subequations}
\begin{align}
\tilde{a}_1^{1,1}(h,h)\triangleq&\frac{1}{2\sigma_x^2}\left[\|h\|^2,\|h\|^2\right]^\t{},\\
\tilde{a}_2^{1,1}(h,h)\triangleq&\frac{1}{2\sigma_x^2}\left[-\|h\|^2,\|h\|^2\right]^\t{},\\
\tilde{a}_3^{1,1}(h,h)\triangleq&-\frac{1}{2\sigma_x^2}\left[\|h\|^2,-\|h\|^2\right]^\t{},\\
\tilde{a}_4^{1,1}(h,h)\triangleq&\frac{3}{2\sigma_x^2}\left[\|h\|^2,\|h\|^2\right]^\t{},
\end{align}
\end{subequations}
and
\begin{subequations}
\begin{align}
\Lambda(h,h)\triangleq&\left[\begin{array}{cc}\|h\|^2& \|h\|^2\\\|h\|^2 &\|h\|^2 \end{array}\right],\\
\overline{\Lambda}(h,h)\triangleq&\left[\begin{array}{cc}\|h\|^2& -\|h\|^2\\-\|h\|^2 &\|h\|^2 \end{array}\right].
\end{align}
\end{subequations}
Using the singularity of $\Lambda(h,h)$ and $\overline{\Lambda}(x,h,h)$, we obtain
\begin{align}
\mu_{x}(1,1,h,h)=&\Ncdf_1\Bigg(0;\frac{\|h\|^2}{2\sigma_x^2},\frac{\|h\|^2}{\sigma_x^2}\Bigg)+\exp\left[\frac{1}{\sigma_x^2}\|h\|^2\right]\Ncdf_1\Bigg(0;\frac{3\|h\|^2}{2\sigma_x^2},\frac{\|h\|^2}{\sigma_x^2}\Bigg).\label{eqn:formulamuxspecial2}
\end{align}
\item \textbf{$\mu_{x}(1,0,h,h)$}: Substituting $s_1=1$, $s_2=0$ and $h_1=h_2=h$ in~\eqref{eqn:formulamux}, we get
\begin{align}
\mu_{x}(1,0,h,h)=&\Ncdf_2\Bigg(\left[\begin{array}{c}0\\0\end{array}\right];\tilde{a}_1^{1,0}(h,h),\frac{1}{\sigma_x^2}\Lambda(h,h)\Bigg)+\Ncdf_2\Bigg(\left[\begin{array}{c}0\\0\end{array}\right];\tilde{a}_2^{1,0}(h,h),\frac{1}{\sigma_x^2}\overline{\Lambda}(h,h)\Bigg)\nonumber\\
&+\Ncdf_2\Bigg(\left[\begin{array}{c}0\\0\end{array}\right];\tilde{a}_3^{1,0}(h,h),\frac{1}{\sigma_x^2}\overline{\Lambda}(h,h)\Bigg)+\Ncdf_2\Bigg(\left[\begin{array}{c}0\\0\end{array}\right];\tilde{a}_4^{1,0}(h,h),\frac{1}{\sigma_x^2}\Lambda(h,h)\Bigg)\label{eqn:formulamuxspecial3}
\end{align}
where
\begin{subequations}
\begin{align}
\tilde{a}_1^{1,0}(h,h)\triangleq&\frac{1}{2\sigma_x^2}\left[\|h\|^2,\|h\|^2\right]^\t{},\\
\tilde{a}_2^{1,0}(h,h)\triangleq&\frac{1}{2\sigma_x^2}\left[-\|h\|^2,\|h\|^2\right]^\t{},\\
\tilde{a}_3^{1,0}(h,h)\triangleq&-\frac{1}{2\sigma_x^2}\left[\|h\|^2,-\|h\|^2\right]^\t{},\\
\tilde{a}_4^{1,0}(h,h)\triangleq&\frac{1}{2\sigma_x^2}\left[\|h\|^2,\|h\|^2\right]^\t{}.
\end{align}
\end{subequations}
Using the singularity of $\Lambda(h,h)$ and $\overline{\Lambda}(x,h,h)$, we obtain
\begin{align}
\mu_{x}(1,0,h,h)=&\Ncdf_1\Bigg(0;\frac{\|h\|^2}{2\sigma_x^2},\frac{\|h\|^2}{\sigma_x^2}\Bigg)+\Ncdf_1\Bigg(0;\frac{\|h\|^2}{2\sigma_x^2},\frac{\|h\|^2}{\sigma_x^2}\Bigg)\\
=&2\Ncdf_1\Bigg(0;\frac{\|h\|^2}{2\sigma_x^2},\frac{\|h\|^2}{\sigma_x^2}\Bigg)\\
=&1-\erf\left(\frac{\|h\|}{2\sqrt{2}\sigma_x}\right).\label{eqn:formulamuxspecial4}
\end{align}

\end{itemize}

\end{document}